\documentclass{article}

\usepackage{caption, amssymb,amsfonts,amsmath, graphicx, geometry}
\usepackage[nice]{nicefrac}
\usepackage{mathrsfs}
\usepackage[colorlinks=true,bookmarksnumbered=true, pdfauthor={Jacob, M\"orters}]{hyperref}

\newcommand{\one}{{\mathsf 1}}

\newcommand{\E}{\mathbb E}
\renewcommand{\P}{\mathbb P}

\renewcommand{\phi}{\varphi}

\newcommand{\N}{\mathbb N}

\renewcommand{\P}{\mathbb P}

\newcommand{\Tex} {T_{\rm ext}}

\newcommand{\U}{\mathcal U}
\newcommand{\Rec}{\mathcal R}
\newcommand{\I}{\mathcal I}
\newcommand{\St}{\mathscr S}
\newcommand{\Co}{\mathscr C}

\newcommand{\heap}[2]  {\genfrac{}{}{0pt}{}{#1}{#2}}
\newcommand{\sfrac}[2] {\mbox{$\frac{#1}{#2}$}}

\newtheorem{theorem}{Theorem}
\newtheorem{lemma}{Lemma}
\newtheorem{proposition}{Proposition}
\newtheorem{remark}{Remark}

\begin{document}
\title{The spread of infections on evolving scale-free networks}

\author{
  Emmanuel Jacob
  \thanks{Ecole Normale Sup\'erieure de Lyon}
  \and
  Peter M\"orters
  \thanks{University of Bath}
}

\maketitle
\begin{abstract}
We study the contact process
on a class of evolving scale-free networks, where each node updates its connections
at independent random times. We give a rigorous mathematical proof that there is a 
transition between a  phase where for all infection rates the infection survives for a long time, at least
exponential in the network size, and a   
phase where for sufficiently  small infection rates extinction 
occurs quickly, at most like the square root of the network size. 
The phase transition occurs when the power-law exponent crosses the value four. This behaviour is in 
contrast  to that of the contact process on the corresponding static model, 
where there is no phase transition, as well as that of a classical mean-field approximation, 
which has a phase  transition at power-law exponent three. 
The new observation behind our result is that temporal variability of networks can simultaneously increase the rate at which
the infection spreads in the network, and decrease the time which the infection 
spends in metastable states.
\end{abstract}

\section{Introduction}

The spread of disease, information or opinion on networks has been the concern of
extensive research over the past decade. Beyond a large body of empirical research and simulation studies
there is now also a growing number of analytic results, based on relatively simple mathematical 
models of the network and the spreading mechanism~\cite{2, 111, 000, 0000}. 
A popular assumption on the networks in which a disease is transmitted
is that they are scale-free, which means that as the network size grows, 
the proportion of nodes  of degree $k$ stabilizes to a limit which, 
as~$k$ increases, 
decays like $k^{-\tau}$ for some positive power-law exponent~$\tau$. Power-law exponents have been statistically estimated 
for a range of networks, and there is an abundance of tractable mathematical models for 
scale-free networks~\cite{01,41,42, 11,12}. 
An established model for the spread of an infection on a network is 
the contact process, or SIS model. In this model every node
can either be infected or healthy.  An infected node
passes the infection at a fixed rate~$\lambda$ to each of its neighbours, and recovers with a fixed rate 
of one. Once it has recovered, it is again susceptible to infection. 

The extinction time of this type of infections on a scale-free network has been subject of a controversial discussion in the literature. 
Pastor-Sattoras and Vespignani~\cite{6,5,8} have used a non-rigorous mean-field calculation to predict a phase transition 
in the behaviour of the contact process. They  predicted that, for networks with a power-law exponent 
$\tau<3$, the infection can survive for a time exponential in the network size regardless of the infection rate.
If $\tau>3$ however, for sufficiently small infection rate, the  time to extinction is at
most polynomial in the network size. Chatterjee and Durrett~\cite{1} have shown that this prediction is 
inaccurate and 
the infection survives for an exponential time for \emph{all} infection rates and power-law exponents. 
An earlier version of this result is due to Berger et al.\ \cite{0}, and a refinement can be found in
Mountford et al.\ \cite{4}. 
The failure of the mean-field approach highlights the need for mathematically
rigorous arguments in this area of complex networks.

Up to now, rigorous mathematical research has focused almost exclusively  
on networks that remain constant while 
the states of the vertices change. In reality, however, connections between individuals change 
over time, and investigating the effect of this temporal variability on processes taking place on 
the network is a question of fundamental importance. 
The study of evolving, or temporal, networks has been identified 
as an important direction for research,  for example in the concluding paragraph of~\cite{2} and in 
the recent survey~\cite{3}. 
The present paper is intended as a first step into territory almost untouched by rigorous mathematics up to now,  offering a rigorous analysis of the contact process on a very simple model of an evolving  scale-free network. 

When setting up our model we make three assumptions, which may be justified, for example, if the network is a human interaction network, the updating models movement or migration of individuals, and infection does not greatly affect mobility of individuals. \emph{First}, we assume that the  time-evolution of the network happens  independently of the infection process. This is 
mainly a simplifying assumption in the light of substantial additional difficulties coming from mutual dependence of processes in random environments.
\footnote{We note however that network dynamics depending on the states of the nodes, are studied under the name adaptive, or co-evolutionary, networks in the nonrigorous  literature, see for example~\cite{21, 20, 22, 23}, and \cite{00} for a rare rigorous example.}
\emph{Second}, we assume that the evolution of the network and the evolution of infection are on the same time scale, which leads to the most interesting interaction between the  evolutions. And \emph{third} we assume that the power of a node, which varies greatly  between nodes in scale-free networks, is not changed over time.

In the absence of individual features of the nodes of the network, there are two natural graph evolutions to study: \emph{Edge-updating} in which the presence or absence of every edge is updated with a fixed rate, and \emph{vertex-updating} in which every vertex refreshes all its connections with a fixed rate. Edge updating has been studied under the name of \emph{dynamical percolation} on a range of regular graph models, see for example~\cite{110, 9}. In our scenario edge-updating would not produce qualitatively different behaviour from the static case, which is why we focus our attention on the 
vertex updating rule. 
In line with our assumptions, vertex updating can serve as a rough model for migration of nodes to 
a different location.
\pagebreak[3]

Our main interest in this paper is in the emergence of qualitatively new phenomena due to the temporal variability of the network. Observe that it is not a priori clear whether temporal variability 
increases or decreases the extinction time of the infection.
Our analysis shows that, on the one hand, the temporal variability is helping the system to exit metastable states, thus \emph{decreasing} the extinction time. As a result, beyond a critical value of $\tau$  a 
phase of fast extinction emerges which is not present for the static model. On the other hand, 
time-variability increases the number of  nodes infected by a single node, leading to a faster spread of the infection, thus 
\emph{increasing} the extinction time. 
As a result, the transition between
the phases of fast and slow extinction 
does not occur at the value $\tau=3$, as in the mean-field calculation~\cite{5}, but at 
the larger value $\tau=4$. This constitutes an interesting new effect emerging
in evolving networks.\\[-7mm]


\section{Statement of the result}

We consider an evolving network model $(G_t)_{t\ge 0}$, a family of graphs indexed by continuous time. Its set of vertices is fixed, for all $t\ge 0$, and 
identified with a set of labels  $V=\{1,\ldots,N\}$. The labels correspond to the strength of vertices, with low labels corresponding to powerful vertices.
At time $t=0$, each unordered pair of vertices $\{x,y\}\subset V$ independently forms an edge with probability \\[-2mm]
$$p_{x,y}= \min\Big\{\frac {\beta N^{2\gamma-1}} { x^\gamma y^\gamma}, 1 \Big\},$$
where $\beta>0$ and $\gamma \in (0,1)$ are the parameters of the model. 
This is a special case of the Chung-Lu model, and it is shown in~\cite{11} that  $G_0$ is a scale-free network with power-law exponent \smash{$\tau=1+\frac1\gamma$}. Moreover, the connection probabilities are asymptotically the same as in all rank-one models including the configuration model, see for example~\cite{41,42, 11}. Note that the expected degree of the vertex labeled~$x$ is proportional to $(N/x)^{\gamma}$,  
reflecting the strong hierarchy of vertices. We remark that, unless $\gamma=\frac12$, the connection probability differs from that in the preferential attachment models, see~\cite{12}, and we believe the results of this paper do not extend to that class.
\pagebreak[3]

The time-evolution of the edges of the network obeys the following rule. Each vertex updates independently with intensity $\kappa>0$.
When vertex~$x$ updates, every unordered pair $\{x,y\}\subset V$, for $y\in V\setminus\{x\}$, forms an edge with probability $p_{x,y}$,
independently of its previous state and of all other edges. The remaining edges $\{w,y\}$ with $w, y\neq x$ remain unchanged. Note that the network evolution is stationary, namely 
for every $t\ge0$, the law of $G_t$ is equal to the law of $G_0$.  We emphasise that the strength of a vertex is given by its label and does not change over time. 

We consider the contact process, or SIS infection, on this evolving network. Every vertex is in either of two states, infected or healthy. 
These states evolve according to the following dynamics.
Every infected vertex turns to healthy with rate~one. 
Every healthy vertex turns to infected 
with rate given by the number of infected neighbours, multiplied by a parameter $\lambda>0$ called the infection rate. The time-evolution of both the network structure and the infection process can be described  
by a continuous time Markov chain. Figure~1 illustrates the principal transitions that affect the state of a given vertex or its set of neighbours.\\[-7mm]

\begin{figure}[h]
\begin{center}
\includegraphics[scale=0.6]{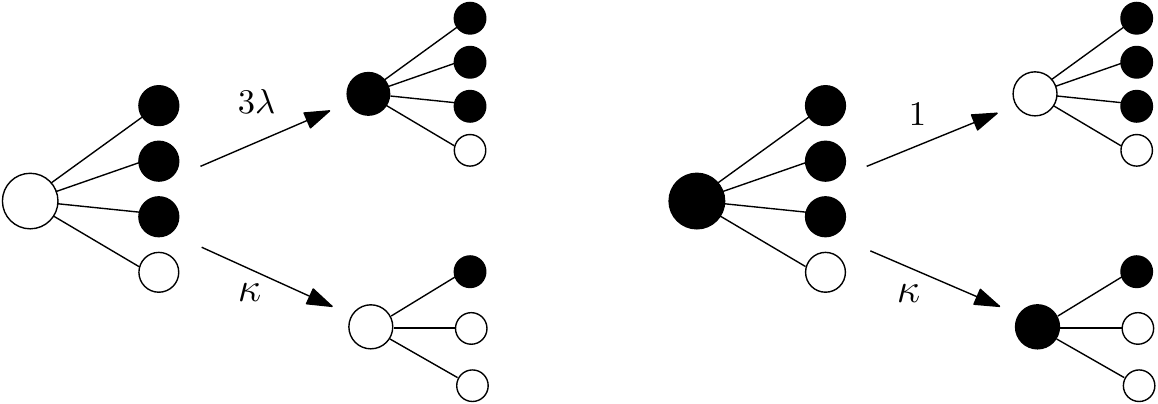}
\end{center}
\caption{Left: The vertex is healthy and has three infected neighbours (represented by black dots) and one healthy one (represented by a white dot, other vertices being not represented). At rate {\footnotesize{$3\lambda$}}, the vertex gets infected , while at rate $\kappa$, it updates, and receives a (possibly) completely new set of neighbours.
Right: The vertex is infected. It turns to healthy at rate 1, and still updates at rate $\kappa$.  
Not represented are updates of other vertices that may lead to loss or creation of new neighbours of the vertex. }
\end{figure}



On a finite graph, no matter how the contact process is started it 
eventually reaches its unique equilibrium, when all vertices are 
healthy and the infection has died. We denote by $\Tex$ the 
(random) extinction time for the infection. The interesting question in this context is now, 
once the infection has broken out, how long it takes until the infection becomes extinct. 
We observe a dichotomy between 
\emph{slow extinction} on the one hand, when $\Tex$ is very large, e.g.\ exponential 
in the number of vertices, and \emph{fast extinction} on the other, when $\Tex$ is small, 
no more than polynomial 
in the number of vertices. It is shown in~\cite{2,4} that
the contact process on the static network $G_0$ exhibits slow extinction for all parameters $\gamma, \beta, \lambda$. Our main result shows that this is different for the evolving network model.

\begin{theorem} \label{MainTheorem}
Consider the contact process on the evolving network $(G_t)_{t\ge 0}$,
and assume that at time $t=0$ every vertex is infected.
\begin{enumerate}
\item[(a)] If $\gamma>1/3$ (or equivalently $\tau<4$), then whatever the other parameters of the network, there exists a positive constant~$c$ such that, uniformly in $N>0$,
$$ \mathbb{P}(\Tex \le e^{cN}) \le e^{-c N}.$$
\item[(b)]  If $\gamma<1/3$ (or equivalently $\tau>4$), then there exists a parameter $\lambda_c>0$ such that, for $\lambda<\lambda_c$, there exists a positive constant $C$ such that, uniformly in $N>0$,
we have
$$ \E[\Tex]\le C \sqrt N.$$
\end{enumerate}
\end{theorem}

\begin{remark}
{\rm  The result of~$(b)$ is true without the assumption on the initial infection. If in~$(a)$ initially not all
vertices are infected, there is additional uncertainty about the outbreak of the infection. In particular, if
under the assumptions of~$(a)$ initially a positive proportion of the vertices are infected, the  event $\Tex > e^{cN}$ still holds with a positive probability, independent of~$N$.}
\end{remark}

We say that a family of events ${\mathscr E}_t$ depending on $t\ge0$ holds \emph{exponentially long}
if there exists $c>0$ such that, for all $N$, 
$$ \mathbb{P}\Big(\bigcap_{s\le e^{c N}}{\mathscr E}_s \Big) \ge 1- e^{-c N}. \\[-1mm]$$
Using this terminology, we can rephrase $(a)$ in Theorem~\ref{MainTheorem} as saying that in
the given regime the infection  survives exponentially long.

\section{Comparison with the static model}

Before entering the rigorous proof of Theorem~\ref{MainTheorem}, we explain, on a heuristic level, how the time-variability of the network can lead to new phenomena that differ both from the rigorous results on static networks, and from the non-rigorous mean-field predictions.\\[-6mm]

\subsection*{Mean-field predictions}
We first describe \emph{informally} the phase-transition as predicted by the mean-field calculations. To demonstrate the difference in behaviour of the static and evolving models, our description focuses on the core assumptions of the predictions and simplifies the actual arguments. Consider the contact process on a static network in which the empirical distribution of the degrees of the vertices is given by some probability measure $\mu$ on $\N= \{1, 2, \ldots\}$. 

The mean-field predictions are based on the principle that instead of taking into account the complex geometry of a network, it suffices to study the behaviour of the contact process around a typical newly infected vertex. The critical infection rate is then the value $\lambda$ at which such a vertex will, on average, transmit the infection once before its recovery. Conditionally on the degree $k$ of the newly infected vertex, the average number of transmissions before recovery is roughly $\lambda k$, therefore to understand criticality it suffices to understand the degree of a newly infected vertex. 

Assuming that the law of the degree of a newly infected vertex is given by $\mu$, one obtains
$$ \lambda_c= \frac 1 {\sum k \mu(k)}.$$
While this assumption is widely used, for example in~\cite{20} in the context of adaptive networks, a possible refinement in our context would take into account that it is easier for a vertex with high degree to become infected. 
%
More precisely, we expect the probability of a vertex being infected to be proportional to the number of possible sources of infection, and hence to its degree. Under this assumption, the law of the degree of a newly infected vertex is the size-biased distribution $\tilde \mu$, defined by
$$ \tilde \mu(k)= \frac {k\mu(k)} {\sum i \mu(i)},$$
which leads naturally to
$$ \lambda_c=\frac 1 {\sum k \tilde \mu(k)} = \frac {\sum k \mu(k)} {\sum k^2 \mu(k)},$$
as in~\cite{8}. If $\mu$ is the asymptotic degree distribution of a scale-free graph with 
power-law exponent~$\tau$, 
we get  $\lambda_c=0$ if $\tau<3$ and $\lambda_c>0$ if $\tau>3$.

The size-biasing effect is crucial in the study of scale-free networks, and can also be found in the rigorous analysis of another infection process, the SIR infection, as in~\cite{0000}.  It often leads to a phase transition
at the value $\tau=3$, see for example~\cite{12,C}, so that the mean-field predictions loosely described here did conform to expectations at their initial appearance 
in~\cite{5}.
As already mentioned, these mean-field predictions have proven to be inaccurate for the contact process on static scale-free networks. It turns out to be essential for the analysis to consider the \emph{geometry} of the network
in the neighbourhood of vertices of large degree, as we explain now.\\[-6mm]

\subsection*{Slow extinction for the static network}
The starting point of the rigorous studies~\cite{1,0,4} concerning the contact process on a static scale-free network is the observation that the infection can survive near a vertex of high degree for much longer than predicted by the mean-field calculations. Specifically, Lemma~(5.3) in~\cite{0} gives an estimate for the survival probability of the contact process on a star graph, which contains only one vertex of degree $k$, and its $k$ neighbours, of degree one. It shows that if $\lambda^2 k$ is larger than an absolute constant, then the time to extinction on this star graph is \emph{exponential} in $\lambda^2 k$. This 
is due to  \emph{quick reinfections} of the central vertex of degree $k$. More precisely, when the central vertex recovers, it is typically surrounded by many infected vertices that {reinfect} it almost instantly. This quick reinfection is not picked up by mean-field calculations, as they keep no memory  of the actual neighbours of the central vertex. \\[-6mm]

\subsection*{Heuristic for the evolving network}
In our evolving network, we can first control how long the \emph{quick reinfections} can maintain the contact process around a vertex, then make \emph{mean-field like} calculations rigorous.
Consider an infected  vertex with large degree $k$. The key observation is that in the evolving network quick reinfections alone can maintain the infection around this vertex for a length of time proportional to the degree~$k$. In fact, the period of quick reinfection will come to an end if either the vertex updates and then recovers before it has infected its new neighbours, or if the vertex recovers and then updates before its old neighbours have reinfected the vertex. Hence, the  infection can be kept in a vertex of degree~$k$ without external reinfection for $\lambda^2 k$ units of time in the evolving network, which has to be compared with the value one (for mean-field calculations not taking reinfections into account) or $\exp(\lambda^2 k)$ (for the static network). 
During the quick reinfection period the vertex will have infected roughly $\lambda^3 k^2$ vertices. 
After the quick reinfection period the vertex can be infected again, but not by a quick reinfection. It will have to be a \emph{new infection} by a vertex that has become a neighbour before infecting it. For such a new infection, we can rigorously argue that the degree of the newly infected vertex follows the size-biased distribution.\pagebreak[3]

Therefore, if $\mu$ stands for the asymptotic degree distribution of the scale-free network, we should observe fast extinction roughly if 
$$ \sum \lambda^3 k^2 \tilde \mu(k)= \lambda^3 \frac {\sum k^3 \mu(k)}{\sum k \mu(k)}<1.$$
Thus, we expect a phase of fast extinction to emerge as soon as $\sum k^3 \mu(k)<\infty$, which corresponds to $\tau>4$ and is in line with \emph{(b)} in Theorem~\ref{MainTheorem}.

\medskip
This heuristic can be made rigorous and leads formally to Lemma~\ref{persist} in the proof of Theorem~\ref{MainTheorem}\emph{\,(a)}. In the proof of Theorem~\ref{MainTheorem}\emph{\,(b)}, the heuristic does not appear explicitly, but has inspired the coupling with a process in which we modify the recoveries, in such a way that a vertex recovers only when we can ensure that it cannot be quickly reinfected by a neighbour.

\section{Graphical representation}

In this section, we provide an equivalent description of our model by a convenient \emph{graphical representation}.
The evolving network model $(G_t)_{t\ge 0}$ is represented 
with the help of the following independent random variables;
\begin{enumerate}
\item[(1)] For each $x\in V$, a Poisson point process of intensity $\kappa>0$, written $\mathcal U^x=(U^x_n)_{n\ge 1}$, 
describing the updating times of~$x$. We also write $\mathcal U^{x,y}=(U^{x,y}_n)_{n\ge 1}$ for the union $\mathcal U^x \cup \mathcal U^y$, 
which is a Poisson point process of intensity $2\kappa$, describing the updating times of the potential edge $\{x,y\}\subset V$. 
\item[(2)] For each potential edge $\{x,y\}\subset V$, a sequence of independent random variables $(C^{x,y}_n)_{n\ge0}$, all Bernoulli with parameter $p_{x,y}$, 
describing the presence/absence of the edge in the network after the successive updating times of the potential edge $\{x,y\}$. More precisely, if $t\ge 0$ then 
$\{x,y\}$ is an edge at time~$t$ if and only if  $C^{x,y}_n=1$ for $n= \vert[0,t]\cap \U^{x,y}\vert$.
\end{enumerate}
Given the network we represent the infection by means of the following set of
independent random variables;
\begin{enumerate}
\item[(3)] For each $x\in V$, a Poisson point process of intensity one, written $\mathcal R^x=(R^x_n)_{n\ge 1}$, describing the recovery times of~$x$.
\item[(4)] For each $\{x,y\}\subset V$, a Poisson point process $\mathcal I_0^{x,y}$ with intensity $\lambda$ describing the infection  times along the edge $\{x,y\}$. Only the trace $\mathcal I^{x,y}$ of this process on the set
$$\bigcup_{n=0}^\infty \{[U^{x,y}_n, U^{x,y}_{n+1}) \colon  C^{x,y}_n=1\} \subset [0,\infty)$$
can actually cause infections. Write $(I^{x,y}_n)_{n\ge1}$ for the ordered points of $\mathcal I^{x,y}$.
If at time $I^{x,y}_n$ vertex $x$ is infected and $y$ is healthy, then $x$ infects $y$. If $y$ is infected and $x$
healthy, then $y$ infects $x$.  Otherwise, nothing happens.
\end{enumerate}
The infection is now described by a process $(X_t(x))_{x\in V}$ with values in \smash{$\{0,1\}^V$}, such that $X_t(x)=1$ if $x$ is 
infected at time~$t$, and $X_t(x)=0$ if $x$ is healthy at time $t$.  
Formally, the infection process $X_t(x)$ associated to this 
graphical representation and to a starting set $A_0$ of infected vertices, is the c\`adl\`ag process with $X_0=\one_{A_0}$ 
evolving only at times $t\in \mathcal R^x\cup \bigcup_{n=1}^\infty I_n^{x,y}$, according to the following rules:
\begin{itemize}
\item If $t\in \mathcal R^x$, then $X_t(x)=0$ (whatever $X_{t-}(x)$). 
\item If $t  \in \mathcal I^{x,y}$, then 
$$
(X_t(x),X_t(y))=
\left\{
  \begin{array}{rl}
    (0,0) & \mbox{ if } (X_{t-}(x),X_{t-}(y))=(0,0). \\
    (1,1) & \mbox{ otherwise.} \\
  \end{array}
\right.
$$
\end{itemize}
 
\section{Proof of Theorem 1(a): Exponential extinction time} 

Fix $\gamma>1/3$ and note that it suffices to look at small values  of
the infection rate~$\lambda>0$. Our strategy is to focus on highly connected vertices, 
which we call \emph{stars}. Informally, when a star is infected typically a proportion of its neighbours is infected at any time. When the star 
recovers, it is likely to be quickly reinfected by its neighbours, keeping the infection alive for a long time. Only when the star updates
before it is reinfected by its neighbours, the recovery is sustained. We use a coupling to look at infections between stars, taking into account
only the sustained recoveries.

More precisely, we partition the vertex set $V$ into two sets
$$ \St:=\{1,\ldots, \lfloor \lambda^{(2+\alpha)/\gamma} N\rfloor\}, \, \Co:=\{ \lfloor \lambda^{(2+\alpha)/\gamma} N\rfloor+1,\ldots, N\},$$
where~$\alpha>0$ is a constant. We call the vertices in $\St$ \emph{stars}, those in $\Co$ \emph{connectors}. Note that a star has an average 
degree of order $\lambda^{-2-\alpha}$, which is large when $\lambda$ is small. The constant $\alpha>0$ is chosen large enough so that we have
\begin{equation}\label{alpha}
\left(\frac 1 \gamma - 3\right)\alpha +  \frac 2 \gamma - 3<0,
\end{equation}
which is possible since $\gamma> 1/3$. The reason for this choice will become clear only at the end of the proof. 
If $x$ and $y$ are stars, then $$p_{x,y}\ge \tilde p_{x,y}:= \beta \lambda^{-4-2\alpha} N^{-1},$$ while if $x$ is a star and $y$ a connector, then 
$$p_{x,y}\ge \tilde p_{x,y}:=\beta\lambda^{-2-\alpha} N^{-1}.$$ By a simple coupling argument, it is enough to prove $(a)$ on the modified network 
where the connection probabilities are replaced by $\tilde p_{x,y}$ and there are no connections between connectors. We 
continue to denote this process by $(X_t(x))_{x\in V}$.

\subsection{Connector update and recovery times} 
In order for the infection to survive for a long time around a star, it requires connectors to infect, and get reinfected by. 
The first lemma ensures that, at any time, 
there are sufficiently many connectors that do not update or recover too quickly.
For any $t\in\N\cup\{0\}$  let
$$
\Co_t:=\{y\in \Co \colon (\Rec^y\cup \U^y)\cap[t,t+2]=\emptyset\}.
$$

\begin{lemma} \label{availableconnectors}
For any sufficiently small $\eta>0$ the event ${\mathscr E}_t:=\{ \vert \Co_t \vert \ge \eta N\}$ holds exponentially long.
\end{lemma}

This lemma is easy and its proof is omitted. From now on, we work  conditionally on the realization of the recovery and update 
times of the connectors, $(\Rec^y)_{y\in \Co}$ and $(\U^y)_{y\in \Co}$, and assume that ${\mathscr E}_t$ holds exponentially long, 
without mentioning explicitly the implied constants. 

\subsection{Extinction time for one star}\
Let $T_\lambda:=\lambda^{-\alpha'}$, where $\alpha'$ is positive constant strictly smaller than $\alpha$, such that
\begin{equation}\label{alpha'}
\left(\frac 1 \gamma - 2\right)\alpha-\alpha' +  \frac 2 \gamma - 3<0,
\end{equation}
which is possible by~\eqref{alpha}. We shall see that a star, alone with its neighbouring connectors, is likely to sustain the infection 
up to time of order $T_\lambda$. 
Note that when $\lambda$ is small, $T_\lambda$ becomes larger, but $\lambda^{-2} T_\lambda$ stays negligible compared to the average degree of a star.

For each star~$x$, let $X^x$ be the infection process obtained by keeping only the edges linking $x$ and $\Co$, with initial condition
$$ 
X_0^x(y)= \left\{
\begin{array}{ll}
0 &\text{if }y\neq x,\\
X_0(x) & \text{if } y=x.
\end{array}
\right.
$$
We call $X^x$ the \emph{infection process around $x$}.
This process depends only on the updating times and recovery times at $x$, and the infection times on edges $\{x,y\}$ for $y\in \Co$. 
Therefore $X^x$ and $X^w$ are conditionally independent if $x$ and $w$ are two different stars, given the realization of  
$(\Rec^y)_{y\in \Co}$ and $(\U^y)_{y\in \Co}$. Moreover, the infection process $X$ satisfies
$$X_t(x) \ge  X_t^x(x) \quad \mbox{for all } t\ge0, x\in \St. $$
We say the infection \emph{persists on $[0,T]$} at a star $x\in \St$, if 
\begin{itemize}
\item $X_T^x(x)=1$, and 
\item $
 \big\vert  \{t\in \U^x\cap[0,T-1] \colon    X_t^x(x)=1,
(\U^x\cup \Rec^x)\cap[t,t+1]=\{t\} \} \big\vert
 \ge \sfrac12 \, \kappa e^{-(1+\kappa)}  \, T.$
\end{itemize}
The second condition ensures that there are enough updating times after which $x$ stays infected for a unit of time without any update or recovery. 
This will be useful later to control interinfections of stars.

\begin{lemma} \label{persist}
Let $x\in\St$ be a star and  ${\mathscr P}_x$ be the event that the infection persists on $[0,T_\lambda]$ at $x$. Then 
$$ \lim_{\lambda \downarrow 0}\liminf_{N\to \infty} \P\big({\mathscr P}_x \,\big|\, X_0(x)=1 \big)=1.$$
\end{lemma}

We say an event depending on $\lambda$ and $N$ holds \emph{with high probability}, if its probability is arbitrarily close to one 
when first $\lambda$ is chosen small enough, and then $N$ sufficiently large. In this sense, Lemma~\ref{persist}, which we prove
in the remainder of this section, states that on $[0,T_\lambda]$ persistence holds with high probability.
A first simple observation is that with high probability we have $\vert \U^x \cap [0,T_\lambda]\vert \le 2 \kappa T_\lambda$ and
\begin{align*}
\big\vert \{t\in \U^x  \cap[0,T_\lambda-1] \colon  (\U^x\cup \Rec^x)  \cap[t,t+1]=\{t\}\} \big\vert 
\ge \sfrac12 \, \kappa e^{-(1+\kappa)} \,T_\lambda.
\end{align*}
It remains to prove that we also have, with high probability, that
$X_{U_n^x}^x(x)=1$ for $n\leq2 \kappa T_\lambda$.
Set 
$U'_0=0$ and recursively \smash{$U'_{n+1}= \min\{1+U'_{n}, \U^x\cap (U'_n,\infty)\}$,} 
so that all time intervals $[U'_n,U'_{n+1}]$ have length at most one. It is enough to prove that, 
with high probability, we have \smash{$X_{U'_n}^x(x)=1$} for $n\leq(1+ 2 \kappa) T_\lambda$.
For $t\le T_\lambda$, define
$$\Co_t(x):=\big\{y\in \Co_{\lfloor t \rfloor} \colon  \{x,y\}\in G_t \big\},$$
the set of neighbours of $x$ that are in $\Co_{\lfloor t \rfloor}$. These are neighbours that we will use to help maintain the infection in $x$.
The following lemma shows that, with high probability, there are enough of them at any time in $[0,T_\lambda]$. 

\begin{lemma} \label{available_neighbours}
$$ \lim_{\lambda\downarrow 0} \liminf_{N\to \infty} \P\big(\vert \Co_{U'_n}(x) \vert \ge \sfrac12 \eta \beta \lambda^{-2-\alpha}\,\, \forall n\le \lfloor (1+ 2 \kappa) T_\lambda \rfloor \big) = 1.$$
\end{lemma}
The result of Lemma~\ref{available_neighbours} follows if  \smash{$\P\{\vert \Co_{U'_n}(x)\vert< \frac12 \eta \beta \lambda^{-2-\alpha}\}$} is,
for large enough $N$, uniformly in $n$,   bounded by a function that is negligible compared to \smash{$\lambda^{\alpha'}$}.
 To this end, observe that, conditionally on $U'_n$, the cardinality of \smash{$\Co_{U'_n}(x)$} is binomial with parameters \smash{$\vert \Co_{\lfloor U'n\rfloor}\vert$} and $\beta\lambda^{-2-\alpha}N^{-1}$. 
As we work on the event \smash{$\{\vert \Co_{\lfloor U'n\rfloor}\vert\geq\eta N\}$}, the cardinality of \smash{$\Co_{U'_n}(x)$} is asymptotically bounded from below by a Poisson variable with parameter 
$\eta \beta \lambda^{-2-\alpha}$. By a standard large deviation bound, the probability that it is less than  $\frac12 \eta \beta  \lambda^{-2-\alpha}$ decays exponentially fast in 
\smash{$\lambda^{-2-\alpha}$}, completing the proof.
\medskip

Using Lemma~\ref{available_neighbours}, we can work conditionally on the event 
$$\{\vert \Co_{U'_n}(x) \vert \ge  \sfrac12 \eta \beta \lambda^{-2-\alpha} \forall n\le \lfloor (1+ 2 \kappa) T \rfloor\}. $$
Lemma~\ref{persist} is proved once we prove the following lemma.

\begin{lemma}\label{short_survival}
For $n\ge0$, we have
$$\P\big(X_{U'_{n+1}}^x(x)=0 \,\big|\, X_{U'_{n}}^x(x)=1 \big) = o(\lambda^{\alpha'}),$$
where the bound is uniform in $n$.
\end{lemma} 

The probability on the left does not depend on $n$, so we may take $n=0$ and $X_0^x(x)=1$. 
On the event \smash{$\{X_{U'_1}^x(x)=0\}$}, we necessarily have $\Rec^x\cap (0,U'_1)\ne \emptyset.$ 
Define random variables 
$$  R=\inf \Rec^x\cap(0, U'_{1}), \qquad L= U'_{1}-\sup \Rec^x\cap(0, U'_{1}).$$
The strategy is to find a connector $y$ that gets infected by~$x$ during the time interval $[0,R]$, and then reinfects $x$ on $[U'_1-L, U'_1]$. 
We look for such a connector in the set $\Co_0(x)$, as \pagebreak[3] these are connected to $x$ during the whole interval
$[0,U'_1]$, and do not recover during this time. For a given $y\in \Co_0(x)$, the probability that the back-and-forth infection between $x$ and $y$ occurs is
$$ \big( 1- e^{-\lambda R}\big) \big( 1- e^{-\lambda L}\big)\le \sfrac 14 \lambda^{2}LR,$$
where the bound uses the inequality $1-e^{-x}\ge x/2$ for $x\in[0,1]$ and the fact that $L, R < 1$.
Hence, the probability that this fails for all $y\in \Co_0(x)$ is bounded by
$$ \left( 1-  \frac {\lambda^{2}LR}4\right)^{\vert \Co_0(x)\vert} \le \exp\left( \sfrac12\eta \beta \lambda^{-2-\alpha} \log( 1-  \sfrac14 \lambda^{2}LR)\right).$$
On the event that the product $LR$ is not too small, ie. \smash{$LR\ge \lambda^{\alpha''}$}, with $\alpha''=\frac12\,(\alpha+\alpha')$, the bound is \smash{$o(\lambda^{\alpha'})$}, as requested.
Finally, to bound the probability that $LR$ is small, we show that
\begin{equation} \label{productbound}
\P\big(LR\le \lambda^{\alpha''}\big)=O\big(\lambda^{\alpha''} \log(\lambda^{-1})\big). 
\end{equation}
Indeed, Equation~\eqref{productbound} follows by bounding separately the probability of the three events 
\smash{$\{L\le 2 \lambda^{\alpha''}\}$}, \smash{$\{R\le 2 \lambda^{\alpha''}\}$}, and
\smash{$\{L\le 1/2, R\le 1/2, LR\le \lambda^{\alpha''}\}$}. The bound on the third event uses the fact that on 
$\{L\le 1/2, R\le 1/2\}$, the random variables $L$ and $R$ are independent and their laws have bounded density on 
$[0,1/2]$. The bound for the first two events also follows from this, we leave the (easy) details to the reader.

\subsection{Reinfections between stars}\ 
Consider an initial condition where the stars in ${\mathscr I}\subset \St$ are infected, and those in $\St\backslash {\mathscr I}$ are healthy. 
From the preceding subsection, we know that, independently for each infected star, the infection is likely to persist at that star. Denote
by ${\mathscr I}'$ the set of infected stars at which the infection persists on $[0,T_\lambda]$.

In this subsection, we prove a quantitative lemma stating that each healthy star has independently some probability of becoming infected 
by some star $s\in {\mathscr I}'$, and then be infected at time $T_\lambda$. We condition on the infection processes $(X^x)_{x\in{\mathscr I} }$, 
which determine ${\mathscr I}'$. 
For $x\in  \St\backslash {\mathscr I}$, we define
$$ T_x:= \inf \big\{t>0 \colon \exists y\in {\mathscr I}', X_t^y(y)=1, \,t\in \I^{x,y}\big\}.$$
Observe that if $T_x$ is finite, then $x$ is necessarily infected at this time, namely we have $X_{T_x}(x)=1$. 
If $T_x\leq T_\lambda$, we define an infection process similar to the process $X^x$, but on the restricted time interval $[T_x,T_\lambda]$. More precisely, define the process \smash{$(\tilde X^x_t)_{t\in [T_x,T_\lambda] }$} as the infection process obtained by keeping only edges between~$x$ and $\Co$, and with initial condition
$$ 
\tilde X_{T_x}^x(y)= \left\{
\begin{array}{ll}
0 &\text{if }y\neq x,\\
1 & \text{if } y=x.
\end{array}
\right.
$$
We further extend this process to $[0,T_\lambda]$ by letting $\tilde X^x_t\equiv 0$ if $t<T_x$ 
(so the process is well-defined even if $T_x>T_\lambda$). It clearly holds that
$$ \forall x\in \St\backslash \mathscr I, X_{T_\lambda}(x) \ge \tilde X^x_{T_\lambda}(x).$$
Moreover, note the following two facts;
\begin{itemize}
 \item The processes $\tilde X^x$, for $x\in \St\backslash \mathscr I$, are independent (conditionally on~$(X^x)_{x\in\mathscr I}$).
\item Conditionally on $T_x<T_\lambda$, the 
probability that $x$ is infected at time $T_\lambda$ goes to 1 when $\lambda$ goes to 0 and $N$ is large. This is obtained by a slight variation of Lemma~\ref{persist}.
\end{itemize}

Therefore, if $\lambda$ is small and $N$ large enough, every healthy star has independently probability at least
$$ \sfrac 12 \P\big(T_x<T_\lambda \ \big\vert \ (X^y)_{y\in \mathscr I} \big)$$
of being infected at time $T_\lambda$. The following lemma bounds this probability 
from below by a quantity depending on $\vert \mathscr I'\vert$.

\begin{lemma} \label{reinfections}
If $x\in \St\backslash \mathscr I$ is a healthy star, ie. $X_0(x)=0$, then
$$ \P\big(T_x<T_\lambda \, \big\vert\, (X^y)_{y\in \mathscr I} \big) \ge 
 1-\exp\big(-\sfrac14 \beta\kappa e^{-(1+\kappa)}  \sfrac {\vert \mathscr I'\vert \lambda^{-3-2\alpha} T} N\big).$$
\end{lemma}

Condition on $(X^y)_{y\in \mathscr I}$, the infection processes around each $y\in \mathscr I$. We first suppose that the infection persists at only one star called $y$, namely $\mathscr I'=\{y\}$. We further condition on the updating times $\U^x$, and hence
on the updating times of the edge $\{x,y\}$, namely $(U^{x,y}_n)_{n\ge 0}$. For $n\ge0$, write 
$$ t_n:=\int_{[U^{x,y}_n, U^{x,y}_{n+1} \wedge T_\lambda]} \one\{X_t^{y}(y)=1\} \, dt, $$ 
the total time when vertex $y$ is infected between consecutive updating times. The probability that $y$ does not infect $x$ is 
\begin{align*}
\prod_{n\ge 0} \big(1- \tilde p_{x,y} (1-e^{-\lambda t_n})\big) 
&\le \prod_{n\ge 0} \big(1-\lambda \tilde p_{x,y} \sfrac12 \min\{1, t_n\}\big)\\
&\le \exp\big(-\sfrac12 (\tilde p_{x,y} \lambda) \sum_{n\ge 0} \min\{1,t_n\} \big).\\[-8mm]
\end{align*}
Using the lower bound for the sum we get from the fact that the infection persists around $y$, 
we continue the inequality by
\begin{align*}
\phantom{\prod_{n\ge 0} \big(1- \tilde p_{x,y} (1-e^{-\lambda t_n})\big) }
&\le \exp\big(-\sfrac14 \kappa e^{-(1+\kappa)} \lambda \tilde p_{x,y} T_\lambda\big) \\
&\le \exp\big(-\sfrac14 \beta\kappa e^{-(1+\kappa)}  \sfrac {\lambda^{-3-2\alpha} T_\lambda} N\big).
\end{align*}
Note that the events that $y$ does not infect $x$, for $y\in \mathscr I'$, are majorised by events 
that are independent given  the realization of $\U^x$. Hence  the probability 
that no vertex of $\mathscr I'$ 
infects $x$ is bounded by the right-hand-side to the power $\vert \mathscr I'\vert$,
and the result follows.

\subsection{Discrete infection process}

We study the infection process on stars at discrete times $0, T_\lambda, 2T_\lambda, \ldots$
For $n\ge0$, we denote by $\mathscr I_n$ the set of stars $x\in \St$ with $X_{nT_\lambda}(x)=1$, and by 
$\mathscr A_n$ as the set of stars~$x\in \mathscr I_n$ such that  the infection persists on $[nT_\lambda, (n+1)T_\lambda]$
at $x$. We start with $\mathscr I_0=\St$, ie.\ initially all stars are infected. 

By Lemma~\ref{persist}, the cardinality of $\mathscr A_0$ is stochastically larger than a binomial random variable with parameters $\vert \St\vert$ and~$\nicefrac12$, if $\lambda$ is small and $N$ large. More generally, for $n\ge 0$, the cardinality of $\mathscr A_n$ conditionally on $\vert \mathscr A_0\vert, \ldots , \vert \mathscr A_{n-1}\vert $ and on $\vert \mathscr I_0\vert, \ldots,\vert \mathscr I_{n}\vert $ is stochastically larger than a binomial random variable with parameters $\vert \mathscr I_n\vert$ and $\nicefrac12$. By Lemma~\ref{reinfections}, conditionally on $\vert \mathscr A_0\vert, \ldots,\vert \mathscr A_{n}\vert $ and on $\vert \mathscr I_0\vert, \ldots,\vert \mathscr  I_{n}\vert$, the cardinality of $\mathscr  I_{n+1}\backslash\mathscr  A_n$ is stochastically larger than a binomial random variable with parameters $\vert \St\vert - \vert \mathscr I_n \vert$ and  
\begin{equation}\label{probab}\sfrac 12\Big(1- \exp\big(-\sfrac14 \beta\kappa e^{-(1+\kappa)}  
\sfrac {\vert A_n\vert  \lambda^{-3-2\alpha} T_\lambda} N\big)\Big).
\end{equation}

Next, we show inductively  that, with high probability, the event $\mathscr E_n:=\{\vert \mathscr I_n\vert\ge \frac14\,\vert \St\vert\}$ holds exponentially long.  
Let $n\geq 0$ and condition on $\mathscr E_n$. Then, as the number of stars is linear in $N$ and each persists on 
$[nT_\lambda, (n+1)T_\lambda]$ independently with probability at least $\nicefrac12$, we infer that the event 
\smash{$\vert \mathscr A_n\vert\ge \frac25 \,\vert \mathscr I_n\vert \ge \frac1{10}\vert \St\vert$} holds with probability at least 
\smash{$1-e^{-cN}$} for some $c>0$, which depends on $\lambda$ but not on $N$. If this holds, then we have
$$ \frac { \vert \mathscr A_n\vert  \lambda^{-3-2\alpha} T_\lambda} N\le 
\sfrac 1 {10}\lambda^{(\frac 1 \gamma - 2)\alpha-\alpha'+\frac 2 \gamma - 3}.$$
The exponent in $\lambda$ is negative, by our choice of $\alpha$ and $\alpha'$, see~\eqref{alpha'}. 
Therefore, choosing $\lambda$ smaller if necessary, we can ensure that the value of \eqref{probab}
is larger than $\nicefrac 13$, which means that each non-infected star, or vertex in $\St\backslash {\mathscr I}_n$, has independently probability at least $\nicefrac 13$ to get infected. This, together with
$\vert \mathscr A_n\vert\ge \frac25\vert \mathscr I_n\vert$, ensures that the event $\vert\mathscr  I_{n+1}\vert \ge \frac14 \, \vert \St \vert$ holds with probability at least \smash{$1-e^{-cN}$} for some (new) constant $c$. Altogether, we have proved that the infection persists exponentially long.

\section{Proof of Theorem~1(b): Polynomial extinction time}

In this section, we propose a coupling of the contact process on an evolving scale-free network, with a \emph{mean-field infection model}. In the mean-field infection model, an infected vertex recovers only when for the original process we can ensure the recovery is sustained and not subject to a quick reinfection. This is the case when the recovery is immediately preceeded by an update of the vertex, so that it is no longer exposed to its original infected neighbours. This allows us to neglect information on the graph structure around the infected vertex prior to the update. Hence we can set up the mean-field infection model independently of the topology of the evolving network, namely as an infection process on the complete graph with appropriately decreased infection rates. 
Finally, the mean-field infection model can be studied with standard techniques, and we can show that the time to extinction is polynomial when $\gamma<1/3$.

\subsection{Mean-field infection model}
To set up the mean-field infection model, for every $x\in V$, we use the point processes $\mathcal R^x$ and $\U^x$ introduced for graphical representation of the original infection process. In addition, for each potential edge $\{x,y\} \subset V$, we define an independent Poisson point process ${\mathcal J}^{x,y}=({ J}^{x,y}_n)_{n\ge 1}$ with intensity $\lambda p_{x,y}$.
The sets ${\mathcal J}^{x,y}$ describe infection times for this process, but the sets $\mathcal R^x$ and $\U^{x}$ have now a different meaning. They describe how recoveries can happen, but in a more involved way. To describe the infection, we use a Markov
process $(Y_t(x))_{x\in V}$ with values in $\{0,1,2\}^V$, where a vertex $x$ can take three different values, the value $0$ still meaning \emph{healthy}, the value~$2$ meaning \emph{infected}, and the value $1$ meaning infected but \emph{ready} to recover. The process evolves according to the following rules;
\vspace{-1mm}
\begin{enumerate}
\item[(1)] If $t  \in { \mathcal J}^{x,y}$, then 
\vspace{-2mm}
$$
(Y_t(x),Y_t(y))=
\left\{
  \begin{array}{rcl}
    (0,0) && \text{ if } (Y_{t-}(x),Y_{t-}(y))=(0,0), \\
    (2,2) && \text{ otherwise.} \\
  \end{array}
\right.
$$
\item[(2)] If $t\in \mathcal U^x$, then 
\vspace{-2mm}
$$
Y_t(x)=
\left\{
  \begin{array}{rcl}
    1 && \text{ if } Y_{t-}(x)=2, \\
    Y_{t-}(x) && \text{ otherwise.} \\
  \end{array}
\right.
$$
\item[(3)] If $t\in \mathcal R^x$, then 
\vspace{-2mm}
$$
Y_t(x)=
\left\{
  \begin{array}{rcl}
    0 && \text{ if } Y_{t-}(x)=1, \\
    Y_{t-}(x) && \text{ otherwise.} \\
  \end{array}
\right.
$$
\end{enumerate}
In other words, once a vertex $x$ is infected, in order to recover, it has to observe an `updating' time in $\U^x$, then a `recovery' time in $\Rec^x$, and no infection time in-between (in $\bigcup_y \mathcal J^{x,y}$).

\subsection{Coupling}

The next proposition states that, if we fix the initial conditions of the infection processes $X$ and $Y$ so that $X_0\le Y_0$, then there exists a coupling maintaining the ordering $X_t\le Y_t$ for all times.

\begin{proposition} Fix deterministic initial conditions $X_0\le Y_0$.
One can construct, on the same probability space, the evolving network $(G_t)_{t\ge 0}$, the infection process \smash{$(X_t(x))_{x\in V, t\ge0}$} on this network, and the mean-field infection process $(Y_t(x))_{x\in V, t\ge 0}$, such that $X_t(x)\le Y_t(x)$ for all $x\in V$ and~$t\ge0$.
\end{proposition}

We stress that the coupling uses the knowledge of the initial condition $X_0$.
As indicated above, we choose the same $\mathcal U^x$ and $\mathcal R^x$ for both models. 
The coupling will actually be a coupling between 
the processes $\mathcal J^{x,y}$ and the processes $\I^{x,y}$, 
which are determined as before by variables $\I_0^{x,y}$,  $C_n^{x,y}$ and $\mathcal U^{x,y}= \mathcal U^x\cup\mathcal U^y$.
Given the latter, the set 
$ {\mathcal J}^{x,y}$ is defined as a subset of ${\mathcal I}_0^{x,y}$, obtained by percolation of parameter $p_{x,y}$, and hence 
is a Poisson point process of parameter $\lambda p_{x,y}$, as requested. However, this percolation procedure will be dependent on the process $X$ and on the $C^{x,y}_n$, as we now explain. We introduce the following notions.
\vspace{-1mm}
\begin{itemize}
\item An \emph{infection time} is an infection time of the original process, namely a time $t\ge 0$ belonging to some $\I^{x,y}$, or equivalently belonging to some $\I_0^{x,y}$ and such that $C_n^{x,y}=1$ if $t\in (U_n^{x,y},U_{n+1}^{x,y})$. 
\item A \emph{potential infection time} is any time $t$ in some $\I_0^{x,y}$.
\item A \emph{true infection time} is an infection time $t$ satisfying the additional condition 
\vspace{-1.5mm}
$$\max\{X_{t-}(x), X_{t-}(y)\}=1.$$
\item  A \emph{potential true infection time} is  a potential infection time~$t$ satisfying additionally
\vspace{-1.5mm}
$$\max\{X_{t-}(x), X_{t-}(y)\}=1.$$
\end{itemize}
It is clear that only true infection times can play any role in the spread of the infection. Observe also that the notion of true infection times depends monotonically on the infection process up to time $t-$, in the sense that any true infection time for~$X$ would also be a true infection time for a process $X'\ge X$. 
For each $x, y$, and $n$, we denote by $F_n^{x,y}$ the first potential true infection time in 
$(U_n^{x,y},U_{n+1}^{x,y})$ if there is any, 
\begin{align*}
 F_n^{x,y}=\inf\big\{t\in \I_0^{x,y} \cap &  (U_n^{x,y},U_{n+1}^{x,y}) \colon 
 \max\{ X_{t-}(x), X_t-(y)\}=1\big\},
\end{align*}
with the convention $\inf \emptyset=+\infty$. The important observation now is that 
the time $F_n^{x,y}$ and the infection process strictly before that time, on $[0, F_n^{x,y})$, are independent of $C_n^{x,y}$. 
This follows from the fact that, on $(U_n^{x,y},U_{n+1}^{x,y})$, before the first potential true infection, there cannot be any true infection along ${x,y}$, whatever the value of $C_n^{x,y}$, and thus the infection process can be determined independently of $C_n^{x,y}$.

To construct $\mathcal J^{x,y}$ let $t\in \I_0^{x,y} \cap (U_n^{x,y},U_{n+1}^{x,y})$. 
If $t$ is the first potential true infection, ie,\ if $t=F_n^{x,y}$, we let $t\in \mathcal J^{x,y}$ if $C_n^{x,y}=1$ and $t\notin\mathcal J^{x,y}$ if $C_n^{x,y}=0$.
Otherwise, sample an independent Bernoulli $p_{x,y}$ random variable to determine whether 
$t$ is kept in $\mathcal J^{x,y}$ or not.

It is clear that $\mathcal J^{x,y}$ has the required law, ie.\ it is a Poisson point process of intensity $\lambda p_{x,y}$ obtained by percolation of parameter $p_{x,y}$ on $\mathcal I_0^{x,y}$. We now explain why the mean-field process $Y$ has to stay above $X$ if it has started with $Y_0\ge X_0$. To this end, we have to ensure two things:
\begin{itemize}
\item If a vertex $x$ recovers for the process $Y$ at some time $t$, then it also recovers for the process $X$ (if it was infected). 
\item If $t\in \mathcal I^{x,y}\cap (U_n^{x,y},U_{n+1}^{x,y})$ is a true infection time for the process $X$, then we have $(Y_t(x),Y_t(y))\ge (1,1)$.  \\[-5mm]
\end{itemize}
The first item is obvious as $Y$ recovers only at times $t\in \Rec^x$, and $X$ also recovers at these times.
To verify the second item  observe that under our assumption we have~$C_n^{x,y}=1$, and hence on $(U_n^{x,y},U_{n+1}^{x,y})$ potential infections and potential true infections coincide with infections and true infections. \pagebreak[3]

Now we distinguish two cases. 
In the first case, $t=F_n^{x,y}$ is the first true infection time for the process $X$ on $(U_n^{x,y},U_{n+1}^{x,y})$. 
By checking all the potential infection times  $t\in \mathcal I_0^{x,y}$ successively
we ensure that $Y$ is above $X$ up to time $t-$. As $t\in \mathcal J^{x,y}$ by construction, we obtain $(Y_t(x),Y_t(y))=(2,2)\ge(1,1)$. 
In the second case, $t$ is not the first true infection time on $(U_n^{x,y},U_{n+1}^{x,y})$. Then we have $F_n^{x,y}<t$, and at time $F_n^{x,y}\in \mathcal J^{x,y}$, we have set the value of $Y$ in $x$ and $y$ to be~2. These values cannot decrease before the next updating time for $x$ or $y$, namely $U_{n+1}^{x,y}$. We deduce again $Y_t(x)= Y_t(y)=2$, and this finishes the proof.

\subsection{Study of the mean-field infection model}

In this part, we show that the infection dies in polynomial time for the mean-field model. We use a simple supermartingale technique, by introducing a process $M(t)$ defined by
$$ M(t):= \sum_{x=1}^N \one_{\{Y_t(x)=2\}} s_2(x) + \sum_{x=1}^N \one_{\{Y_t(x)=1\}} s_1(x).$$
The values $s_1(x)$ and $s_2(x)$, for $x\in\{1,\ldots, N\}$, are nonnegative scores, associated to infected vertices, which should roughly quantify to which extent the fact that $x$ is infected contributes to the survival of the infection. Actually, we want to define them in such a way that, if the infection rate $\lambda$ is small,  $M$ is a tractable supermartingale with respect to the filtration \smash{$({\mathscr F}_t)$} generated by $(Y_t)$. Informally, $s_2(x)$ should be larger than $s_1(x)$, and, for $i=1, 2$, we should have $s_i(x)$ larger than $s_i(x')$ if $x<x'$, as $p_{x,y}\ge p_{x',y}$ for all $y\ne x,x'$. We have
\begin{align*}
\frac 1 {dt}  \E\big[M(t+dt)-M(t)\big\vert {\mathscr F}_t\big]
&=
 \sum_{x\colon Y_t(x)=2}\kappa (s_1(x)-s_2(x)) \\
&\quad+ \sum_{x \colon Y_t(x)=1} \left( -s_1(x) + \sum_{y \colon y\ne x} \lambda p_{x,y} (s_2(x)-s_1(x))\right)\\
&\quad+ \sum_{x \colon Y_t(x)=0} \left( \sum_{\heap{y \colon y\ne x,}{ Y_t(y)\ge 1}} \lambda p_{x,y} s_2(x)\right).
\end{align*}
In this equality we have written for each vertex $x$, the change of its contribution to the total score, namely $s_1(x)\one{\{Y_t(x)=1\}}+s_2(x)\one{\{Y_t(x)=2\}},$ 
induced by the infinitesimal transition probabilities of the infection process changing the state of~$x$. 
 The third sum can be rewritten as
$$ \sum_{x \colon Y_t(x)\ge 1} \sum_{\heap{y \colon y\ne x,}{ Y_t(y)=0}} \lambda p_{x,y} s_2(y)
\le \sum_{x \colon Y_t(x)\ge 1} \sum_{y} \lambda p_{x,y} s_2(y).$$
Writing $S(x):=\sum_{y} p_{x,y}$ and $T(x):=\sum_{y} p_{x,y} s_2(y),$ we get the following bound,
\begin{eqnarray*}
\frac 1 {dt}  \E\big[M(t+dt)-M(t)\big\vert {\mathscr F}_t\big] 
&\le& \sum_{x \colon Y_t(x)=2}\lambda T(x)- \kappa \big(s_2(x)-s_1(x)\big)\\
&&  +  \sum_{x \colon Y_t(x)=1}\lambda T(x)+\lambda S(x) \big(s_2(x)-s_1(x)\big)- s_1(x).
\end{eqnarray*}

It is easy to bound $S(x)$, which is the average degree of $x$, as 
$$ S(x)\le \frac  \beta {N^{1-2\gamma} x^\gamma}\sum_{y=1}^N y^{-\gamma}\le \frac \beta {1-\gamma} \left(\frac N x\right)^{\gamma}.$$
We now choose $s_1$ and $s_2$ as 
$$ s_1(x)= \left(\frac N x\right)^{2\gamma} \qquad \qquad s_2(x)=s_1(x)+ \left(\frac N x\right)^\gamma.$$
Then, we can bound $T(x)$ as 
$$T(x) \le  \frac {2\beta} {N^{1-4\gamma} x^\gamma}\sum_{y=1}^N y^{-3\gamma}\le \frac {2\beta} {1-3\gamma} \left(\frac N x\right)^\gamma,$$
where we used $\gamma< \nicefrac 13$. Therefore, for any $x$, we have
$$\lambda T(x)- \kappa \big(s_2(x)-s_1(x)\big)\le \left( \frac {2\lambda \beta}{1-3\gamma} - \kappa\right) \left(\frac N x\right)^\gamma$$
and
\begin{align*}
\lambda T(x) & +\lambda S(x) \big(s_2(x)-s_1(x)\big)- s_1(x) \le \frac {2\lambda \beta}{1-3\gamma}\left(\frac N x\right)^\gamma 
 + 
\left( \frac {\lambda \beta}{1-\gamma} - 1\right) \left(\frac N x\right)^{2 \gamma}.
\end{align*}
If $\lambda$ is small enough, the negative terms compensate all the positive ones, and we get, for some positive constant $\nu>0$,
\begin{align*}
\frac 1 {dt} \E\big[M(t+dt)-M(t) \big\vert {\mathscr F}_t \big] 
&\le- \nu \left(\sum_{x \colon Y_t(x)=2} s_2(x)^{\frac12}+ \sum_{x \colon Y_t(x)=1} s_1(x)\right) \\
&\le - \nu \sqrt{M(t)}.
\end{align*}
We introduce the process $Z(t)=\sqrt{M(t)}+\nu t/2$, and get, on the event $M(t)>0$,
\begin{align*}
\frac 1 {dt} \E\big[ & Z(t+dt)  -Z(t)\big\vert {\mathscr F}_t \big] 
 \le \frac \nu 2+\frac 1 {2\sqrt{M(t)}} \, \frac 1 {dt}\E\big[M(t+dt)-M(t)\big\vert {\mathscr F}_t \big]\le 0.
\end{align*}
Observing that the extinction time $T_{\rm ext}$ of the process $Y$ is also the hitting time of 0 for the process $M$, we get that 
$Z(t\wedge T_{\rm ext})$ defines a positive supermartingale, converging to $Z(T_{\rm ext})= \frac\nu2 T_{\rm ext}$. Using the optional stopping theorem we deduce $\E[Z(T_{\rm ext})]\le Z(0)$, and therefore
$$ \E[T_{\rm ext}]\le \frac 2 \nu \sqrt{M(0)}
\leq  \frac 2 \nu \sqrt{ {\textstyle\sum_x s_2(x)}} \leq \frac 2 \nu \sqrt{ \frac{2N}{1-2\gamma}}.$$
Hence the expectation of $T_{\rm ext}$ grows at most like $\sqrt N$
and the second part of Theorem~\ref{MainTheorem} is proved.

\bigskip
\noindent \emph{EJ was supported by CNRS, and PM was supported by EPSRC grant EP/K016075/1.}

\end{document}